\documentclass{gtart}

\def\ifplaintex{\expandafter\ifx\csname documentclass\endcsname\relax}

\ifplaintex 
\else
\fi

\def\gtm{{\mathsurround=0pt\it $\cal G\mskip-2mu$eometry \&\ 
$\cal T\!\!$opology $\cal M\mskip-1mu$onographs}}    

\def\gtp{{\mathsurround=0pt\it $\cal G\mskip-2mu$eometry \&\ 
$\cal T\!\!$opology $\cal P\!$ublications}}  

\def\recd{{\small Received:\qua\receiveddate\ifx\reviseddate\relax
\else\qquad Revised:\qua\reviseddate\fi\par}} 

\def\volumenumber#1{\def\thevolumenumber{#1}}
\def\volumeyear#1{\def\thevolumeyear{#1}}
\def\volumename#1{\def\thevolumename{#1}}
\def\papernumber#1{\def\thepapernumber{#1}}
\def\pagenumbers#1#2{\def\startpage{#1}\def\finishpage{#2}}
\def\published#1{\def\publishdate{#1}}
\def\received#1{\def\receiveddate{#1}}
\def\revised#1{\def\reviseddate{#1}}
\def\accepted#1{\def\accepteddate{#1}}
\def\asciititle#1{\def\theasciititle{#1}}
\def\covertitle#1{\def\thecovertitle{#1}}

\def\asciiaddress#1{\def\theasciiaddress{#1}}

\long\def\asciiabstract#1{\long\def\theasciiabstract{#1}}

\def\shorttitle#1{\def\theshorttitle{#1}}

\let\\\par
\let\thevolumenumber\relax\let\thepapernumber\relax
\let\thevolumeyear\relax\let\startpage\relax
\let\finishpage\relax\let\publishdate\relax\let\receiveddate\relax
\let\reviseddate\relax\let\accepteddate\relax\let\theasciititle\relax
\let\thecovertitle\relax\let\theasciiauthors\relax\let\theasciiaddress\relax
\let\theasciiabstract\relax

\let\theerratum\relax\let\theasciiemail\relax
\let\theshortauthors\relax\let\theshorttitle\relax

\def\startpage{1}\def\finishpage{15}\def\thepapernumber{77}

\volumenumber{2}
\volumename{Proceedings of the Kirbyfest}
\volumeyear{1999}

\long\def\maketitlep{   

\count0=\startpage

\gtm\nl        
{\small Volume \thevolumenumber: \thevolumename\nl 
\ifx\theerratum\relax\else Erratum \erratumnumber\nl\fi
Pages \startpage--\finishpage\nl}

\vglue 0.1truein   

{\parskip=0pt\leftskip 0pt plus 1fil\def\\{\par\smallskip}{\ifplaintex\large
\else\Large\fi\bf\thetitle}\par\medskip}   
\vglue 0.05truein 

{\parskip=0pt\leftskip 0pt plus 1fil\def\\{\par}{\sc\theauthors}
\par\medskip}%
 
\vglue 0.03truein 

{\small\leftskip 25pt\rightskip 25pt{\bf Abstract}\stdspace\theabstract

{\bf AMS Classification}\stdspace\theprimaryclass
\ifx\thesecondaryclass\relax\else; \thesecondaryclass\fi\par
{\bf Keywords}\stdspace \thekeywords\par}\vglue 7pt

}   

\font\phead=cmsl9 scaled 950
\font=cmsl9 scaled 1050
\font\pnum=cmbx10 scaled 913
\font\lnum=cmbx10 
\font\pfoot=cmsl9 scaled 950
\font=cmsl9 scaled 1050
\ifplaintex
\headline{\vbox to 0pt{\vskip -4.5mm\line{\small\phead\ifnum
\count0=\startpage ISSN 1464-8997 (on line)
1464-8989 (printed) \hfill {\pnum\folio}\else\ifodd\count0\def\\{ }%
\ifx\theshorttitle\relax\thetitle\else\theshorttitle\fi\hfill{\pnum\folio}
\else\def\\{ and }{\pnum\folio}\hfill\ifx\theshortauthors\relax\theauthors
\else\theshortauthors\fi\fi\fi}\vss}}
\footline{\vbox to 0pt{\vglue 0mm\line{\small\pfoot\ifnum\count0=\startpage
Published \publishdate:\qua\copyright\ \gtp\hfill\else
\gtm, Volume \thevolumenumber\ (\thevolumeyear)\hfill\fi}\vss
}}
\else
\makeatletter
\def\@oddhead{{\small\lhead\ifnum\count0=\startpage ISSN 1464-8997 (on line)
1464-8989 (printed) \hfill {\lnum\number\count0}\else\ifodd\count0
\def\\{ }\ifx\theshorttitle\relax \thetitle \else\theshorttitle\fi\hfill
{\lnum\number\count0}\else\def\\{ and }{\lnum\number\count0}
\hfill\ifx\theshortauthors\relax 
\theauthors\else\theshortauthors\fi\fi\fi}}\def\@evenhead{@oddhead}
\def\@oddfoot{\small\lfoot\ifnum\count0=\startpage Published \publishdate:\qua\copyright\ \gtp\hfill\else
\gtm, Volume \thevolumenumber\ (\thevolumeyear)\hfill\fi}
\def\@evenfoot{@oddfoot}
\makeatother
\fi

\let\maketitlepage\maketitlep
\let\makeshorttitle\maketitlepage
\let\maketitle\maketitlepage

\newwrite\gtoutfile
\long\gdef\makeheadfile{  
{\def\\{, }\def\s{ }
\immediate\openout\gtoutfile head.xxx
\immediate\write\gtoutfile{To: math@arxiv.org}
\immediate\write\gtoutfile{Subject: put OR rep NNNNN:ppppp}
\immediate\write\gtoutfile{--text follows this line--}
\immediate\write\gtoutfile{Proxy-for: \ifx\theasciiauthors\relax
\theauthors\else\theasciiauthors\fi\s<\ifx\theasciiemail\relax\theemail\else\theasciiemail\fi>}
\immediate\write\gtoutfile{\noexpand\\}
\immediate\write\gtoutfile{Authors: \ifx\theasciiauthors\relax
\theauthors\else\theasciiauthors\fi}
{\def\\{ }\immediate\write\gtoutfile{Title: \ifx\theasciititle\relax
\thetitle\else\theasciititle\fi}}
\immediate\write\gtoutfile{Subj-class: GT or SG, GR etc}
\immediate\write\gtoutfile{MSC-class: \theprimaryclass\ifx\thesecondaryclass\relax\else, \thesecondaryclass\fi}
\immediate\write\gtoutfile{Journal-ref: Geom. Topol. Monogr. \thevolumenumber\s
(\thevolumeyear) \startpage-\finishpage}
\immediate\write\gtoutfile{Comments: Published by Geometry and Topology Monographs at}
\immediate\write\gtoutfile{\s\s\s  http://www.maths.warwick.ac.uk/gt/GTMon\thevolumenumber/paper\thepapernumber.abs.html}
\immediate\write\gtoutfile{\noexpand\\}
\immediate\write\gtoutfile{}
\ifx\theasciiabstract\relax
\immediate\write\gtoutfile{\theabstract}\else
\immediate\write\gtoutfile{\theasciiabstract}\fi
\immediate\write\gtoutfile{}
\immediate\write\gtoutfile{\noexpand\\}
\immediate\write\gtoutfile{}
\immediate\closeout\gtoutfile}}  

\def\maketitlepage{\maketitlep\makeheadfile}
\let\makeshorttitle\maketitlepage
\let\maketitle\maketitlepage

\volumenumber{4}
\volumename{Invariants of knots and 3-manifolds (Kyoto 2001)}
\volumeyear{2002}
\papernumber{18}  
\pagenumbers{281}{294}
\received{30 November 2001}
\revised{7 March 2002}
\accepted{22 July 2002}
\published{13 October 2002}

\usepackage{amssymb,amsmath,amscd} 
\usepackage{graphicx}
\usepackage{floatflt}

\newtheorem{thm}{Theorem}[section]  
        
\newtheorem{cor}[thm]{Corollary}
\newtheorem{prop}[thm]{Proposition}
\theoremstyle{definition}
\newtheorem{defn}[thm]{Definition}  
          
\newtheorem{exm}[thm]{Example}

\begin{document}
\title{$(2+1)$-dimensional topological quantum field 
theory\\\vglue-7pt\\with 
a Verlinde basis and Turaev-Viro-Ocneanu\\invariants of $3$-manifolds}
\covertitle{$(2+1)$-dimensional topological quantum field 
theory\\with 
a Verlinde basis and Turaev-Viro-Ocneanu\\invariants of $3$-manifolds}
\asciititle{(2+1)-dimensional topological quantum field theory with 
a Verlinde basis and Turaev-Viro-Ocneanu invariants of 3-manifolds}
\shorttitle{TQFT with Verlinde basis and Turaev-Viro-Ocneanu TQFT }

\authors{Nobuya Sato\\Michihisa Wakui}                  
\address{Department of Mathematics and Information Sciences\\Osaka 
Prefecture University, Sakai, Osaka, 599-8531, 
Japan\\\smallskip\\Department of Mathematics, Osaka University\\Toyonaka, 
Osaka 560-0043, Japan}

\asciiaddress{Department of Mathematics and Information Sciences\\Osaka
Prefecture University, Sakai, Osaka, 599-8531, Japan\\Department of
Mathematics, Osaka University\\Toyonaka, Osaka 560-0043, Japan}

\email{nobuya@mi.cias.osakafu-u.ac.jp, wakui@math.sci.osaka-u.ac.jp}

\begin{abstract}   
In this article, we discuss a $(2+1)$-dimensional topological 
quantum field theory, for short TQFT,  with a Verlinde basis. As a 
conclusion of this general theory, we have a Dehn surgery formula.
We show that Turaev-Viro-Ocneanu TQFT has a Verlinde basis. Several 
applications of this theorem are exposed. Based on Izumi's data of 
subfactors, we list several computations of Turaev-Viro-Ocneanu 
invariants for some 3-manifolds. 
\end{abstract}

\asciiabstract{
In this article, we discuss a (2+1)-dimensional topological 
quantum field theory, for short TQFT,  with a Verlinde basis. As a 
conclusion of this general theory, we have a Dehn surgery formula.
We show that Turaev-Viro-Ocneanu TQFT has a Verlinde basis. Several 
applications of this theorem are exposed. Based on Izumi's data of 
subfactors, we list several computations of Turaev-Viro-Ocneanu 
invariants for some 3-manifolds.}


\maketitle

\section{Introduction}
At the beginning of the 1990's, Turaev-Viro-Ocneanu TQFT \cite{Ocneanu} was 
introduced by A. Ocneanu by using a type II$_1$ subfactor with finite index 
and finite depth as a generalization of Turaev-Viro TQFT \cite{TuraevViro} 
which was derived from the quantum group $U_q(sl(2,\Bbb{C}))$ at certain roots 
of unity. K. Suzuki and the second author \cite{SuzukiWakui} have found a 
Verlinde basis for the Turaev-Viro-Ocneanu TQFT from an $E_6$-subfactor and 
have computed the invariant for several $3$-manifolds including lens spaces 
$L(p,q)$, where $p, q$ are less than or equal to $12$. An interesting 
consequence of their research is that the Turaev-Viro-Ocneanu invariant 
distinguishes orientations for specific manifolds. It is known that the 
Turaev-Viro invariant cannot distinguish orientations since for an orientable 
closed 3-manifold it coincides with the square of absolute value of the 
Reshetikhin-Turaev invariant \cite{Turaev}. These facts let us expect that 
Turaev-Viro-Ocneanu invariant has  a lot of more useful information than the 
Turaev-Viro invariant.

After this introduction, two expository sections of subfactors and the 
Turaev-Viro-Ocneanu invariant are in order. 

In Section 4, we introduce the notion of a $(2+1)$-dimensional TQFT with 
a \lq\lq Verlinde basis''. This is a special case of a $(2+1)$-dimensional 
TQFT $Z$ with extra assumptions for the representations of the action of the 
mapping class group of $S^1 \times S^1$ on $Z(S^1 \times S^1)$. From this 
general theory, we can afford a useful Dehn surgery formula of 3-manifold 
invariants. Actually, we have many useful formulas for calculating the 
Turaev-Viro-Ocneanu invariant of some concrete $3$-manifolds, for example, 
lens spaces and Brieskorn 3-manifolds. 

In Section 5, we state a theorem that the Turaev-Viro-Ocneanu TQFT has a 
Verlinde basis in our sense. Hence, we have the Dehn surgery formula. We list 
several applications derived from this theorem.

In Section 6, we show that $S$ and $T$-matrices in Izumi's sector theory 
associated with the Longo-Rehren subfactor and $S$ and $T$-matrices in 
the Turaev-Viro-Ocneanu $(2+1)$-dimensional TQFT are complex conjugate. 
Namely, for any subfactor there is a complex conjugate isomorphism as matrix 
representations of $SL_2(\Bbb{Z})$. 
In the last part of this article, we concretely compute 
Turaev-Viro-Ocneanu invariants of lens spaces and homology $3$-spheres based 
on Izumi's data \cite{Izumi2} of $S$ and $T$-matrices.

Finally, we would like to thank Y. Kawahigashi and M. Izumi for valuable 
discussions and comments. The first author is partially supported by 
Grants-in-Aid for Scientific Research, Encouragement of Young Scientists, by 
JSPS.  

\section{Subfactors and fusion algebras}
In this section, we will have a quick exposition of subfactor 
theory. For more precise treatment, see \cite{EK} for instance.
Let us start with the definition of von Neumann algebras. Let $H$ be 
a Hilbert space and $B(H)$ be the set of bounded linear 
operators on $H$. Recall that $B(H)$ is closed under the adjoint operation:
for an operator $a \in B(H)$, there exists the operator $a^*$ such that 
$\langle a \xi, \eta \rangle= \langle \xi, a^* \eta \rangle $ holds for every 
vectors $\xi$, $\eta$ $\in H$. With this involution, $B(H)$ becomes a 
$*$-algebra. Let $M$ be a $*$-subalgebra of $B(H)$. This $M$ is a {\it von 
Neumann algebra} if $M$ is closed under the weak topology. (Namely, 
in the weak topology,  a sequence of operators $\{a_n \}$ converges to $a$ if 
for arbitrary vectors $\xi$, $\eta$ in $H$, $\langle a_n \xi, \eta \rangle 
\to \langle a \xi, \eta \rangle$  holds. ) 
A von Neumann algebra $M$ is called a {\it factor} if $M \cap M'$ is 
$\Bbb{C} 1_M$, where $1_M$ is the identity element of $M$ and $M'$ means 
the set of operators which commute with each element of $M$ in $B(H)$. 
Remark that every von Neumann algebra is decomposed into a direct integral
of factors. Factors are classified to three types, type I, type II and 
type III. We are mainly interested in type II and type III. (von Neumann 
factors belonging to Type I  are of the form  $B(H)$ for some $H$.) We do 
not get into further on factors, but we mainly use a type II$_1$ factor, 
which is a type II factor with  the unique normal normalized trace. 

As we see in the following example, we can have {\it subfactors} of a factor. 
\begin{exm}
Let $A_n$ be $\bigotimes_{k=1}^n M_2(\Bbb{C})$ and $B_n$ be 
$ \left( 
\begin{array}{cc}
 1 & 0 \\
 0 & 1
\end{array}
\right) 
\otimes
\bigotimes_{k=1}^{n-1} M_2(\Bbb{C})$. 
Then trivially $B_n$ is a subalgebra of $A_n$. Besides, we can embed 
$A_n$ into $A_{n+1}$ sending $x$ to $x \otimes 1$. In a similar way, 
embed $B_n$ into $B_{n+1}$. 
Then, taking the inductive limit of $A_n$ and taking the weak closure 
induced from the trace on $\cup_{n=1}^\infty A_n$, we have a factor of 
type II$_1$ $A=\overline{\cup_{n=1}^\infty A_n}$ and a subfactor of $A$, 
$B=\overline{\cup_{n=1}^\infty B_n}$. This is a simple example of 
subfactors. Those factors are approximated by the finite dimensional 
algebras. In such a case,  we call our factor AFD(= Approximately 
Finite Dimensional) in brief. 
\end{exm}

Now, let us define Jones index for a type II$_1$ subfactor $N \subset M$. 
The bigger factor can be considered as a left $N$-module ${}_N M$ and this 
left module is projective. So it can be projectively decomposed. Assume 
the number of direct summands is finite. (Otherwise we think that Jones 
index takes the value infinity.) Then, left $N$-module ${}_N M $ 
is isomorphic to ${}_N N \oplus \cdots \oplus {}_N N \oplus {}_N Np $, 
where $p$ is a projection in $N$. Then, Jones index $[M:N]$ is defined 
to be (the number of ${}_N N$ in the direct sum) + tr$(p)$, where tr is 
the unique trace on $M$. For instance, $[A:B]=4$ in the previous example.  
Jones index is sometimes denoted by $\mbox{\rm dim} {}_N M$. 

For a subfactor $N \subset M$, note that $M$ can be viewed as $N$-$M$,  
$M$-$N$, $M$-$M$ and $N$-$N$ bimodule. Using this fact, we are going to 
produce new bimodules. Let ${}_P X_Q$ be a $P$-$Q$ bimodule, where $P, 
Q \in \{ M, N \}$. Then, we define the $Q$-$P$ bimodule ${}_Q \bar{X}_P$ 
in the following manner. As a Hilbert space, it is the conjugate 
Hilbert space of $X$. The $Q$-$P$ action is defined by $q \cdot 
\bar{\xi} \cdot p=\overline{p^* \cdot \xi \cdot q^*}$, where $\xi \in X$, 
$p \in P$ and $q \in Q$. This $Q$-$P$ bimodule ${}_Q \bar{X}_P$ is called 
the conjugate bimodule of ${}_P X_Q$. For instance, $\overline{{}_N M_M} 
= {}_M M_N$. Another operation is the relative tensor product of bimodules. 
This is similar to the tensor product of algebraic bimodules over a certain 
ring except the completion as a Hilbert space. 
We say that a $P$-$Q$ bimodule ${}_P X_Q$ is irreducible if the set of 
bounded  $P$-$Q$ linear maps $\mbox{End}({}_P X_Q)$ is isomorphic to 
the scalar $\Bbb C$. Set $g={}_N M_M$. Take relative tensor products 
$g \underset{M}{\otimes} \bar{g} \underset{N}{\otimes} \cdots 
\underset{N}{\otimes} g$. Then this can be decomposed into the irreducible 
bimodules $\bigoplus_i m_i \;  {}_N {X_i}_M $, where $m_i$ is the multiplicity 
of the irreducible bimodule ${}_N {X_i}_M $ in the relative tensor products. 
(Recall that this is very much similar to representation theory 
of a compact group.) Finally, the dimension function $[{}_P X_Q]$ of 
${}_P X_Q$ is defined to be $[{}_P X_Q] =(\mbox{dim} {}_P X)(\mbox{dim} 
X_Q) \in [1, \infty]$. (Note that this corresponds to the square of 
`` quantum dimension'' in different  literatures.)

A set of four kinds of irreducible bimodules with the dimension function 
obtained from a subfactor $N \subset M$ with finite Jones index, closed under 
the operations, relative tensor product, conjugation and direct sum is called 
the {\it graded fusion rule algebra} associated with $N \subset M$. 
If the number of unitary equivalence classes of four kinds of irreducible 
bimodules is finite, then subfactor $N \subset M$ is said to have {\it finite 
depth}. 

For type III factors, it is more natural to use endomorphisms than bimodules. 
Let $M$ be a type III factor. Take endomorphisms $\rho$, $\sigma$ $\in 
\mbox{End}(M)$. We have a natural unitary equivalence relation between  
two endomorphisms in such a way that $u\rho(x)u^*=\sigma(x)$, where 
$u$ is a unitary in $M$. Equivalence classes of $\mbox{End}(M)$ are 
called {\it  sectors}. For sectors there are notions 
of direct sum, composition($=$ \lq\lq tensor product''), conjugation, 
irreducibility and statistical dimension($=$ \lq\lq quantum dimension''), 
which play the roles of the corresponding operations 
for the bimodules obtained from a type II$_1$ subfactor. In both cases of type 
II$_1$ subfactors and type III factors, it is a certain $C^*$-tensor category 
that we have finally constructed.

\section{Turaev-Viro-Ocneanu $(2+1)$-dimensional TQFT}

Let $N \subset M$ be a subfactor of type II$_1$ with finite Jones index 
and finite depth. As we have seen in the previous section, we have the 
graded fusion rule algebra associated with this subfactor. Then, a {\it 
quantum 6j-symbol} is defined as a composition $\xi_4 \cdot (\xi_3 \otimes id) 
\cdot (id \otimes \xi_1)^* \cdot \xi_2^*$ of the following four 
intertwiners, i.e., bounded $M$-$M$ linear homomorphisms 
\[
\begin{CD}
X \underset{M}{\otimes} A \underset{M}{\otimes} Y 
@> id \otimes \xi_1>>X \underset{M}{\otimes}B \\
@V{\xi_3 \otimes id}VV  @VV{\xi_2}V \\
 C \underset{M}{\otimes} Y 
@> \xi_4 >>D
 \end{CD}
\]
where, $X$, $Y$, $A$, $B$, $C$ and $D$ are irreducible $M$-$M$ bimodules 
and $\xi_i's$ are orthonormal basis of intertwiners. (For irreducible 
bimodules $X$, $Y$ and $Z$, we introduce an inner product $\langle \xi, 
\eta \rangle $ for two intertwiners $\xi$ and $\eta$ in 
$\mbox{Hom}(X \underset{M}{\otimes} Y, Z)$ by $\xi \cdot \eta^*$.)

Let us consider a quantum $6j$-symbol as a value of a tetrahedron in the 
following way. In the picture, each vertex, edge and face correspond 
to II$_1$ factor $M$, an irreducible bimodule and an intertwiner, 
respectively. 

\vspace {-10pt}
\cl{\setlength{\unitlength}{0.7cm}
\mbox{}\hspace{50pt}\begin{picture}(21,6)
\put(1.0,1.9){\scalebox{1.1}[1.1]{\includegraphics[width=3cm]{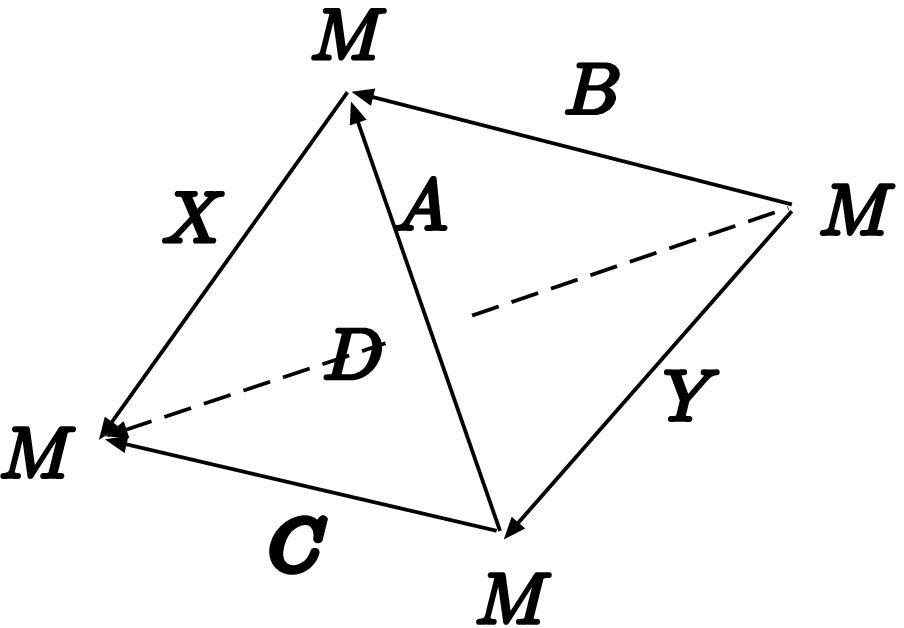}}}
\put(5.3,3){$=\ [B]^{-\frac{1}{4}}[C]^{-\frac{1}{4}}
\ \xi _4\cdot (\xi _3\otimes id)\cdot
(id\otimes \xi _1)^{\ast }\cdot \xi _2^{\ast }$}
\end{picture}}
\vspace {-35pt}

Let $V$ be a closed 3-manifold with triangulation $\cal T$ and let $\cal M$ 
be the $C^*$-tensor category of $M$-$M$ bimodules obtained from a subfactor 
$N \subset M$. Then, the Turaev-Viro-Ocneanu invariant $Z_{\cal M}(V, 
{\cal T})$ is defined in the following way. First, assign a II$_1$ factor $M$ 
for each vertex. Then, assign an irreducible bimodule for each edge. 
Finally, assign an intertwiner for each face of a tetrahedron in $\cal T$. 
We denote the set of edges in $\cal T$ by $E$. We denote all the 
possible assignment of edge colorings and face colorings by $e$ and 
$\varphi$, respectively.  Note that we have a complex number for each 
tetrahedron which is induced from a quantum $6j$-symbol. We denote this 
value of a tetrahedron $\tau$ by $W(\tau;e,\varphi)$. Multiply all tetrahedra 
and take summations over all possible assignments of face colorings and edge 
colorings:
$$
Z_{\cal M}(V,{\cal T})=\lambda^{-a} \sum_e(\prod_E [X]^{\frac{1}{2}} )
\sum_\varphi \prod_\tau W(\tau;e,\varphi)
$$
where $[X]=(\mbox{dim}{}_M X)(\mbox{dim} X_M)$, $\lambda=\sum_{{}_M X_M} 
[X]$(called the {\it global index}) and $a$ is the number of vertices in 
$\cal T$.

This $Z_{\cal{M}}(V,\cal{T})$ turns out to be independent of the 
choice of triangulations. Hence, it is a topological invariant of $V$ 
called the {\it Turaev-Viro-Ocneanu invariant}. We may simply write 
$Z_{\cal{M}}(V)$ for this invariant. One can extend this construction to 
a (2+1)-dimensional unitary TQFT. We again use $Z_{\cal M}$ to express the 
Turaev-Viro-Ocneanu TQFT. 

\section{$(2+1)$-dimensional TQFT with a Verlinde basis}
\par 
In this section we briefly explain a general theory of $(2+1)$-dimensional TQFT with a Verlinde basis. For a $(2+1)$-dimensional TQFT with a Verlinde basis, we have a Dehn surgery formula of the $3$-manifold invariant induced from the TQFT. By applying the Dehn surgery formula to some specific $3$-manifolds including lens spaces and Brieskorn 3-manifolds, we get formulas such that the invariants of them can be calculated from the $S$ and $T$-matrices associated with the TQFT. 
\par 
Let $Z$ be a $(2+1)$-dimensional TQFT and $\{ v_i\} _{i=0}^m$ a basis of $Z(S^1\times S^1)$. Then a family of framed link invariants $J(L;i_1,\cdots ,i_r)\in \mathbb{C},\ i_1,\cdots ,i_r\in \{ 0,1,\cdots ,m\} $ is defined as follows.
\par 
Let $L$ be a framed link in $S^3$ with $r$-components $L_1,\cdots ,L_r$ and $X$ the link exterior of $L$ in $S^3$, that is $X=\overline{S^3-N(L_1)\cup \cdots \cup N(L_r)}$, where $N(L_i)$ is a tubular neighbourhood of $L_i$. 
Then we have a cobordism $W_L:=(X;\partial N(L_1)\cup \cdots \cup \partial N(L_r),\emptyset )$. This cobordism induces a linear map $Z_L:Z(S^1\times S^1)^{\otimes m}\longrightarrow \mathbb{C}$. By using this linear map $Z_L$, a complex number $J(L;i_1,\cdots ,i_r)$ is defined by 
$$J(L;i_1,\cdots ,i_r):=Z_L(v_{i_1}\otimes \cdots \otimes v_{i_r})$$
for each $i_1,\cdots ,i_r\in \{ 0,1,\cdots ,m\} $. It is easily seen that $J(L;i_1,\cdots ,i_r)$ is an invariant of framed links in $S^3$. 

\begin{floatingfigure}[r]{4.0cm}
\cl{\includegraphics[height=4cm]{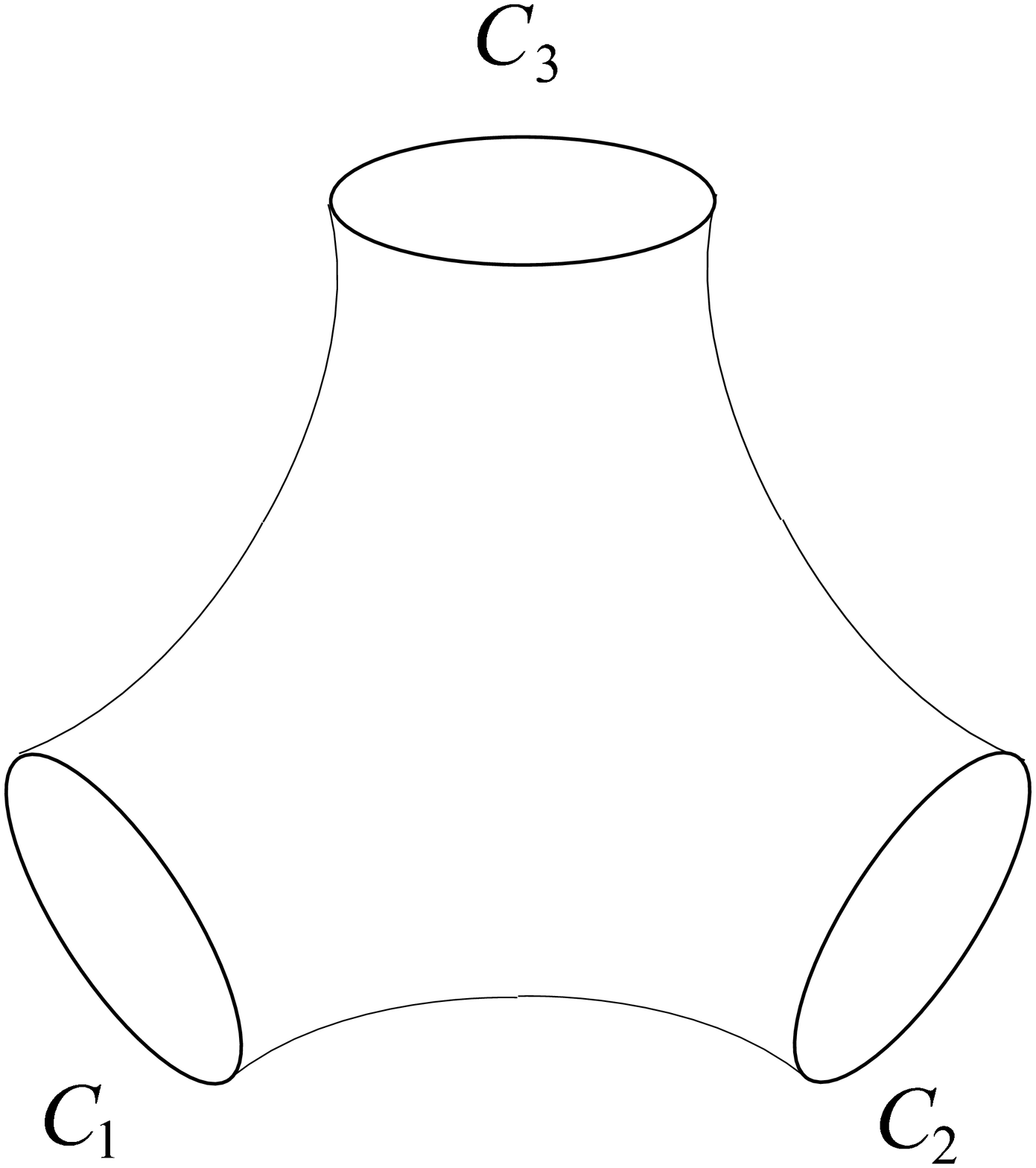}}
\end{floatingfigure}
Let $Z$ be a $(2+1)$-dimensional TQFT. Taking the cobordism $W_{\text{fusion}}:=(Y\times S^1; \Sigma _1\sqcup \Sigma _2, \Sigma_3)$, where $Y$ is the compact oriented surface 
depicted in the figure and $\Sigma _i=C_i\times S^1$ for $i=1,2,3$, 
we get a linear map $Z_{W_{\text{fusion}}}: Z(S^1\times S^1)\otimes Z(S^1\times S^1)\longrightarrow
Z(S^1\times S^1)$ from the axiom for TQFT \cite{Atiyah, Turaev}. This map $Z_{W_{\text{fusion}}}$ defines an associative and unital product on $Z(S^1\times S^1)$. It can be easily shown that this algebra is commutative and the identity element of it is given by
$Z_{W_0}(1)$, where $W_0:=(D^2\times S^1; \phi , S^1\times S^1)$. We call this algebra the {\it fusion algebra associated with} $Z$.

Any $(2+1)$-dimensional TQFT $Z$ induces an action of the mapping class group of an oriented closed surface $\Sigma $ on the vector space $Z(\Sigma )$ \cite{Atiyah, Turaev}. In case of $\Sigma =S^1\times S^1$, the mapping class group of it is isomorphic to the group 
$SL_2(\Bbb{Z})$ of integral $2\times 2$-matrices with determinant $1$, that is generated by $S=\begin{pmatrix} 0 & -1 \\ 1 & 0
\end{pmatrix}$ and $T=\begin{pmatrix} 1 & 1 \\ 0 & 1\end{pmatrix}$ with relations
$S^4=I,\ (ST)^3=S^2$. The matrices $S$ and $T$ correspond to the orientation preserving homeomorphisms from $S^1\times S^1$ to $S^1\times S^1$ which are defined by $S(z,w)=(w,\bar{z})$ and $T(z,w)=(z,zw),\ (z,w)\in S^1\times S^1$, respectively, where we regard $S^1$ as the set of complex numbers of absolute value $1$. 
The actions of $S$ and $T$ on $Z(S^1\times S^1)$ are denoted by $Z(S)$ and $Z(T)$, respectively.

\begin{defn}
Let $Z$ be a unitary $(2+1)$-dimensional TQFT. 
An orthonormal basis $\{ v_i\} _{i=0}^m$ of $Z(S^1\times S^1)$ is called a {\it Verlinde basis} if the following conditions are satisfied.
\begin{enumerate}
\item[(i)]
$v_0$ is the identity element of the fusion algebra. 
\item[(ii)]
$Z(S)$ is presented by a unitary and symmetric matrix, and $Z(T)$ is presented by a unitary and diagonal matrix with respect to $\{ v_i\} _{i=0}^m$. 
\item[(iii)]
$Z(S)^2v_i\in \{ v_i\} _{i=0}^m$ for all $i\in \{ 0,1,\cdots ,m\} $, and $Z(S)^2v_0=v_0$. 
\item[(iv)]
We write $Z(S)v_i=\sum\limits_{j=0}^mS_{ji}v_j,\ S_{ji}\in \mathbb{C}$ for all $i\in \{ 0,1,\cdots ,m\} $. Then 
\item[\rlap{\hglue -10pt(a)}] $S_{0i}\not= 0$ for all $i\in \{ 0,1,\cdots ,m\} $
\item[\rlap{\hglue -10pt(b)}] $N_{ij}^k:=\sum\limits_{l=0}^m\dfrac{S_{il}S_{jl}\overline{S_{lk}}}{S_{0l}}$ is a non-negative integer for all $i,j,k \in \{ 0,1,\cdots ,m\} $
\item[\rlap{\hglue -10pt(c)}] $\{ N_{ij}^k\} _{i,j,k=0,1,\cdots ,m}$ coincide with the structure constants of the fusion algebra with respect to $\{ v_i\} _{i=0}^m$ : $v_iv_j=\sum\limits_{k=0}^mN_{ij}^kv_k.$
\end{enumerate}
\end{defn}

By using the gluing axiom \cite{Atiyah, Turaev} in TQFT, we have the following proposition. 

\begin{prop}
Let $Z$ be a $(2+1)$-dimensional TQFT with Verlinde basis $\{ v_i\} _{i=0}^m$. 
Let $M$ be a closed oriented $3$-manifold. If $M$ is obtained from $S^3$ by Dehn surgery along a framed link $L=L_1\cup \cdots \cup L_r\subset  S^3$, then $Z(M)$ is given by the formula
$$Z(M)=\sum\limits_{i_1,\cdots ,i_r=0}^mS_{i_1,0}\cdots S_{i_r,0}J(L;i_1,\cdots ,i_r).$$
\end{prop}

We apply the above proposition to the case where $M$ is the lens space 
$L(p,q)$ , $q=1,2$ or the Brieskorn 3-manifold 
$M(p,q,r)=\{ (u,v,w)\in \mathbb{C}^3\ \vert \ u^p+v^q+w^r=0,\ {|u|}^2+{|v|}^2+{|w|}^2=1\} $, where $p,q,r\geq 2$. 
Since these manifolds are obtained from $S^3$ by Dehn surgery along framed links depicted in the figure on the next page, we obtain the following formulas. 
\par \quad 
$Z(L(p,1))=\sum\limits_{i=0}^mt_i^{p}S_{i0}^2\quad (p\in \mathbb{N})$,
\par \quad 
$Z(L(p,2))=\sum\limits_{i,j=0}^mt_i^{\frac{p+1}{2}}t_j^2S_{i0}S_{j0}S_{ij}$\quad $(p\in \mathbb{N}$ is odd),
\par \quad 
$Z(M(p,q,r))=\sum\limits_{i,j,k,l=0}^mt_i^pt_j^qt_k^rt_l\dfrac{S_{i0}S_{j0}S_{k0}S_{il}S_{jl}S_{kl}}{S_{l0}}$
\par\noindent  
where $Z(T)v_i=t_iv_i\ (i=0,1,\cdots ,m)$. 

\par \medskip 
These formulas are helpful to compute the values of Turaev-Viro-Ocneanu 
invariants from several concrete subfactors in the section 6.

\hbox{}\quad\includegraphics[height=3.2cm]{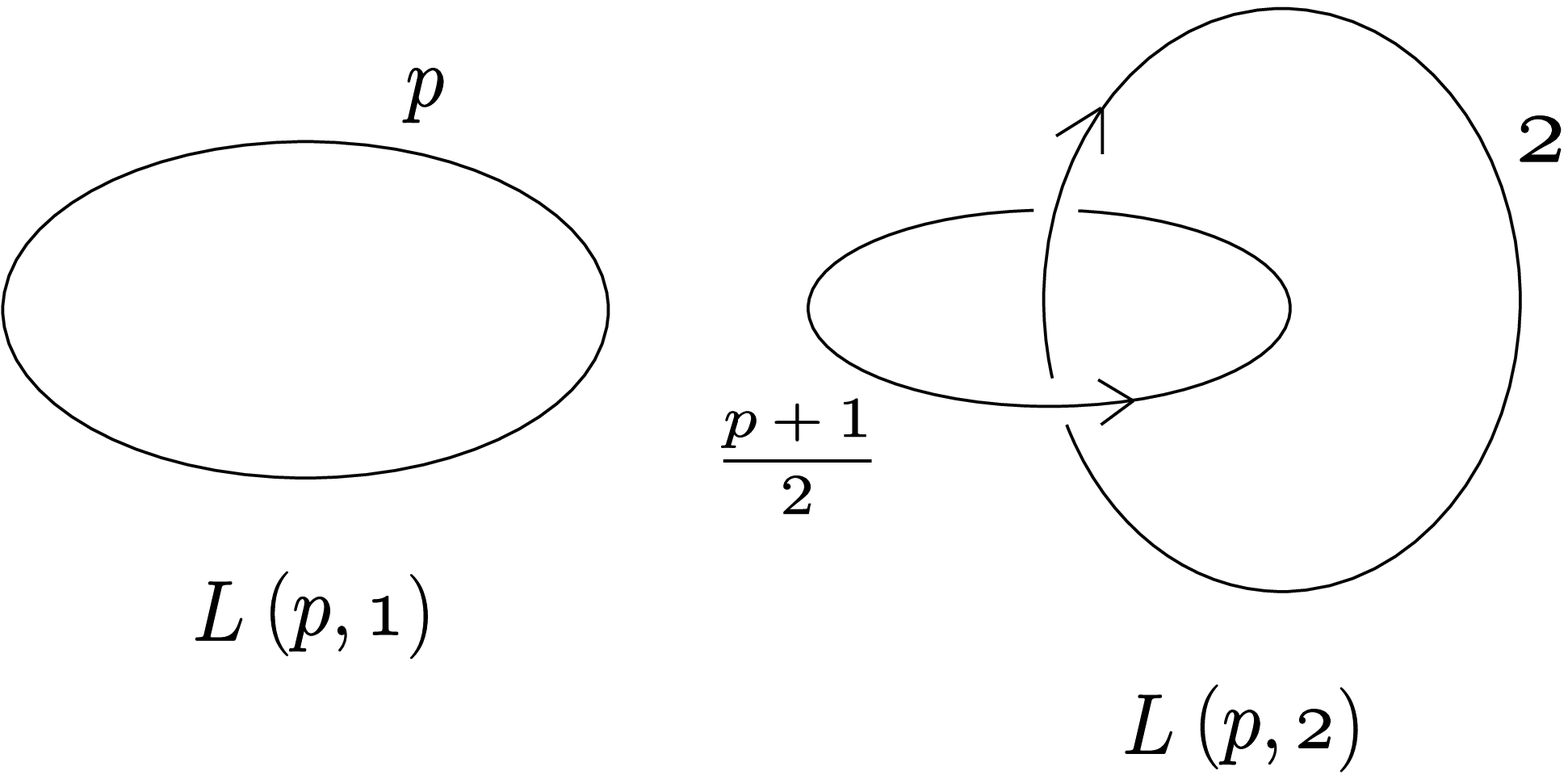}\qquad\quad
\raise 20pt \hbox{\includegraphics[height=3.2cm]{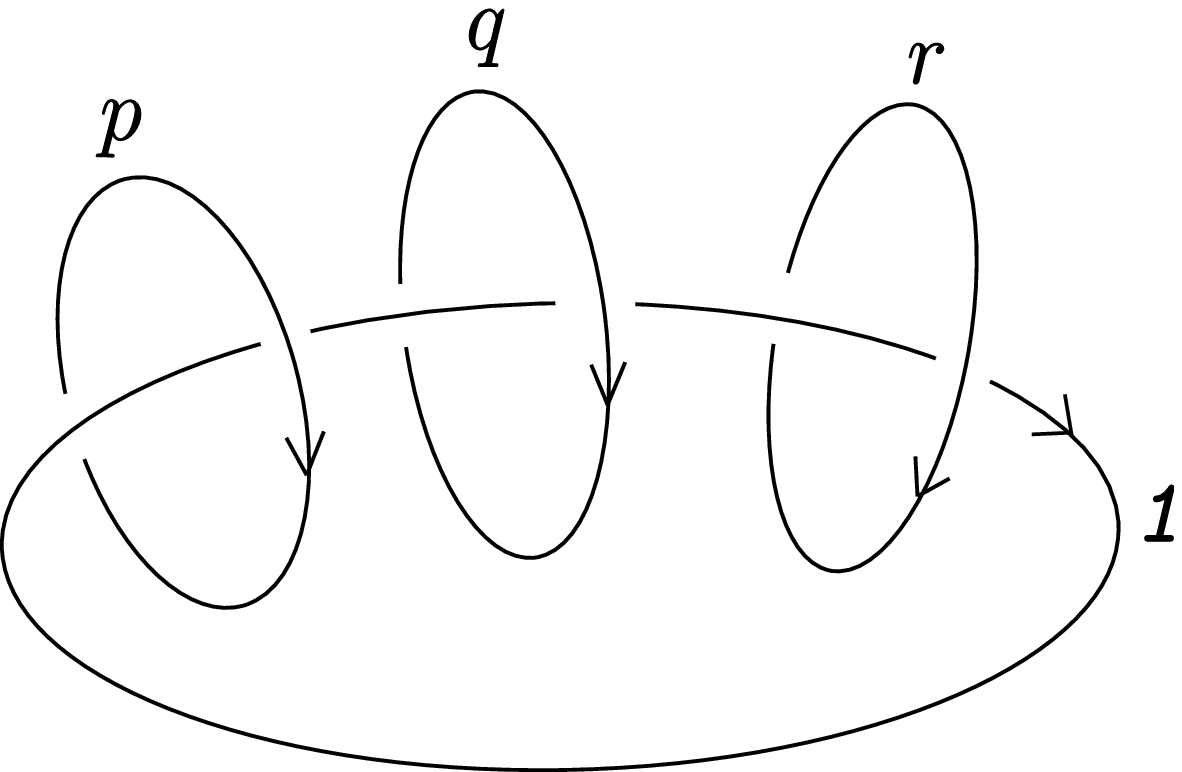}}

\vglue -10pt
{\small \hbox{\hbox{}\hspace{3.7in} $M(p,q,r)$}}

\section{Applications to Turaev-Viro-Ocneanu TQFT}
Let $N \subset M$ be an AFD II$_1$ subfactor with finite Jones index 
and with finite depth. Then, we can construct a new subfactor out of 
this, which is called the {\it asymptotic subfactor} 
$$M \vee M^{op} \subset \vee_{n=0}^\infty M_n=M_\infty$$
where $N \subset M \subset M_1 \subset \cdots \subset M_n \subset 
\cdots$ is the Jones tower, $\vee$ means \lq\lq generated by'' and 
$M^{op}$ is the opposite algebra of $M$.
Let us denote the $C^*$-tensor category of $M_\infty$-$M_\infty$ 
bimodules obtained from the asymptotic subfactor by ${\cal M}_\infty$. 
It turns out that ${\cal M}_\infty$ is  modular in the sense of Turaev. 
It is sometimes easier to treat another tensor category than 
${\cal M}_\infty$ itself. Instead of ${\cal M}_\infty$, we introduce the 
{\it tube system} associated with $N \subset M$, which is equivalent to 
${\cal M}_\infty$ as modular categories. Originally, the tube system 
was introduced by  A. Ocneanu to analyze the asymptotic subfactors 
\cite{Ocneanu}, \cite{O}.

First, we introduce the {\it tube algebra} $\mbox{Tube}({\cal M})$, 
due to Ocneanu.
As a linear space $\mbox{Tube}({\cal M})$ is equal to $\oplus_{X,Y,A} 
\mbox{Hom} (X \underset{M}{\otimes} A,A \underset{M}{\otimes} Y)$. The 
product structure is given by Turaev-Viro-Ocneanu invariant of the join 
$A\ast S^1$ of the annulus $A$ and the circle $S^1$ as follows:

\vglue 10pt
\cl{\small
\setlength{\unitlength}{0.635cm}
\begin{picture}(21.2,3.5)
\put(1,1){\framebox(2,2)}
\put(1.8,0.4){$A$}
\put(1.8,3.2){$A$}
\put(0.38,2.0){$X$}
\put(3.2,2.0){$Y$}
\put(1.7,0.85){$\ll $}
\put(1.7,2.87){$\ll $}
\put(0.8,2){$\vee $}
\put(2.82,2){$\vee $}
\put(1.8,2){$\xi $}
\put(3.6,1.9){$\cdot$}
\put(4.8,1){\framebox(2,2)}
\put(5.6,0.4){$B$}
\put(5.6,3.2){$B$}
\put(3.85,2.0){$X '$}
\put(7.0,2.0){$Z$}
\put(5.5,0.85){$\ll $}
\put(5.5,2.87){$\ll $}
\put(4.6,2){$\vee $}
\put(6.62,2){$\vee $}
\put(5.6,2){$\eta $}
\put(7.45,1.9){$=\delta _{Y,X '}
\frac{[Y]^{\frac{1}{2}}}{[X]^{\frac{1}{4}}[Z]^{\frac{1}{4}}}$}
\put(11.1,1.9){$\displaystyle\sum\limits_{\xi ,\eta ,\zeta }$}
\put(12.1,1.9){$Z\Biggl ($}
\put(13.2,0.1){\includegraphics[width=2.25cm]{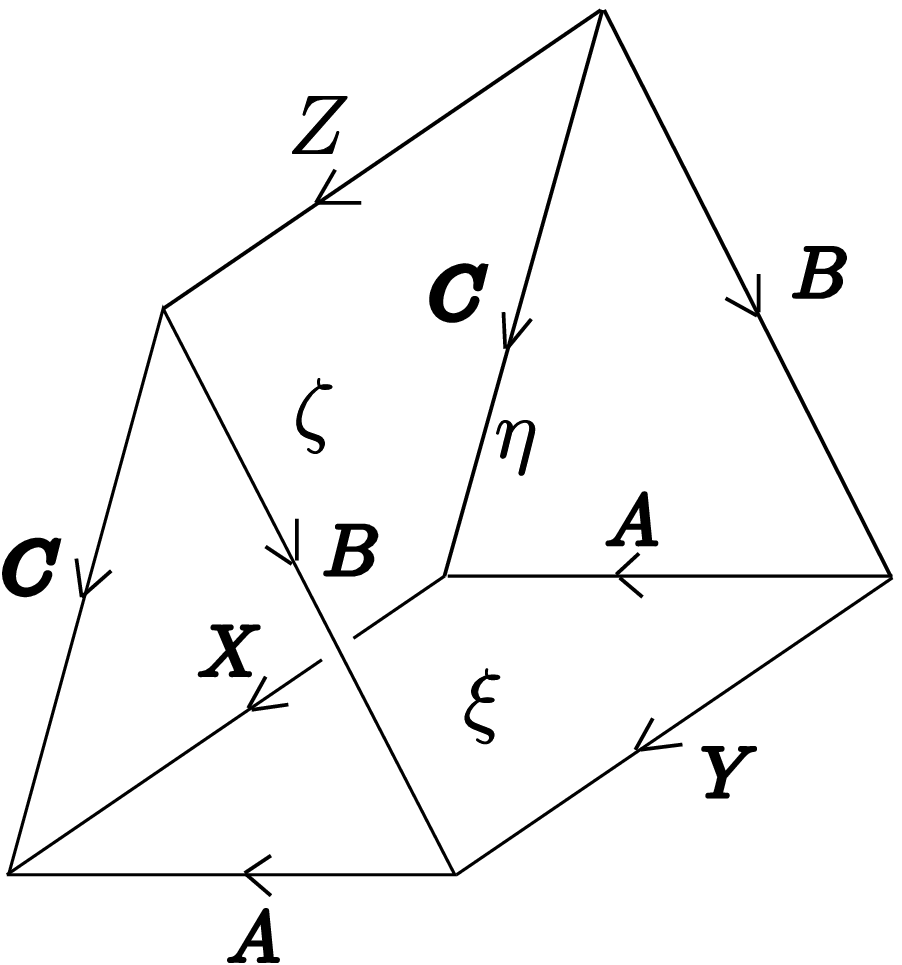}}
\put(17,1.9){$\Biggr )$}
\put(18,1){\framebox(2,2)}
\put(18.8,0.4){$C$}
\put(18.8,3.2){$C$}
\put(17.35,2.0){$X$}
\put(20.2,2.0){$Z$}
\put(18.7,0.85){$\ll $}
\put(18.7,2.87){$\ll $}
\put(17.8,2){$\vee $}
\put(19.82,2){$\vee $}
\put(18.8,2){$\zeta $}
\end{picture}}

\vglue -10pt

Indeed, it turns out  that $\mbox{Tube}({\cal M})$ is a finite 
dimensional $C^*$-algebra. So, $\mbox{Tube}({\cal M}) \cong 
\oplus_{i=0}^d M_{n_i}(\mathbb{C})$. Let $\{ p_0, \cdots, p_d \}$ be 
minimal projections in each direct summand of $\mbox{Tube}({\cal M})$. 

\begin{defn}
 The {\it tube system} associated with $N \subset M$ is defined 
to be a modular category with objects, morphisms and braidings in the 
following: \\
Objects are the $\mathbb{C}$-linear span of minimal projections of 
Tube$({\cal M})$. \\
Morphisms are defined in the following way:\\
Hom$(p_i, p_j)=
\left\{ 
\begin{array}{cr}
 \Bbb{C} & i=j \\
 \{0 \}  & i \ne j 
\end{array}
\right.$

An element of Hom$(p_i \otimes p_j, p_k)$ is a vector in the 
Hilbert space of a surface in the following picture. 
$$\includegraphics[width=3.3cm]{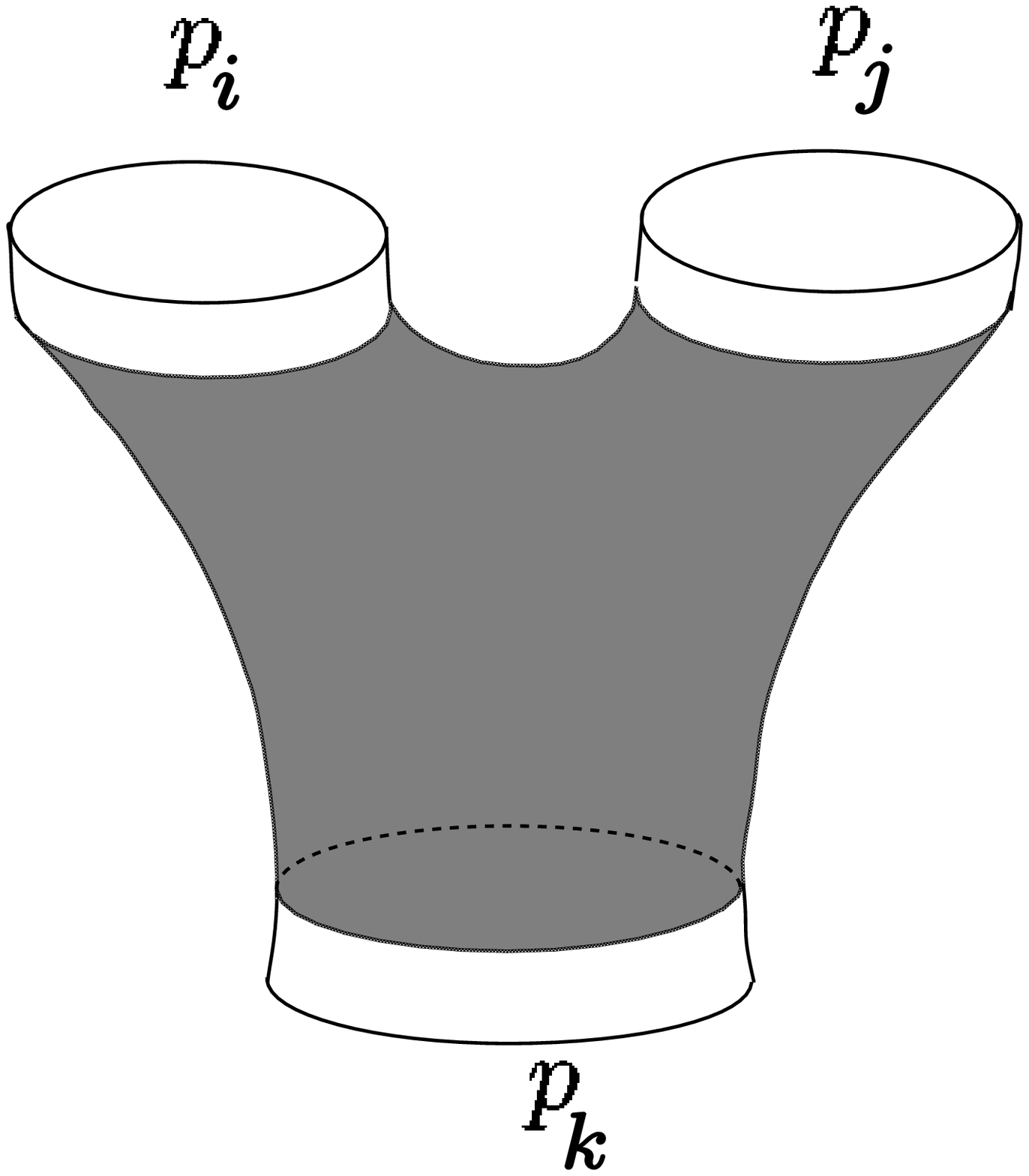}$$
Finally, a unitary braiding $\epsilon_{i,j} \in \mbox{\rm Hom}
(p_i \otimes p_j,p_j \otimes p_i )$ is defined by the 
following formula:

$$\setlength{\unitlength}{0.65cm}
\begin{picture}(21,3.5)
\put(4.0,2.0){$\epsilon _{ij}\ =\ \displaystyle \sum\limits_{k ,\xi ,\eta }$}
\put(7,2){$Z_{\mathcal{M}}\Biggl ($}
\put(8.7,0.5){\includegraphics[width=2.6cm]{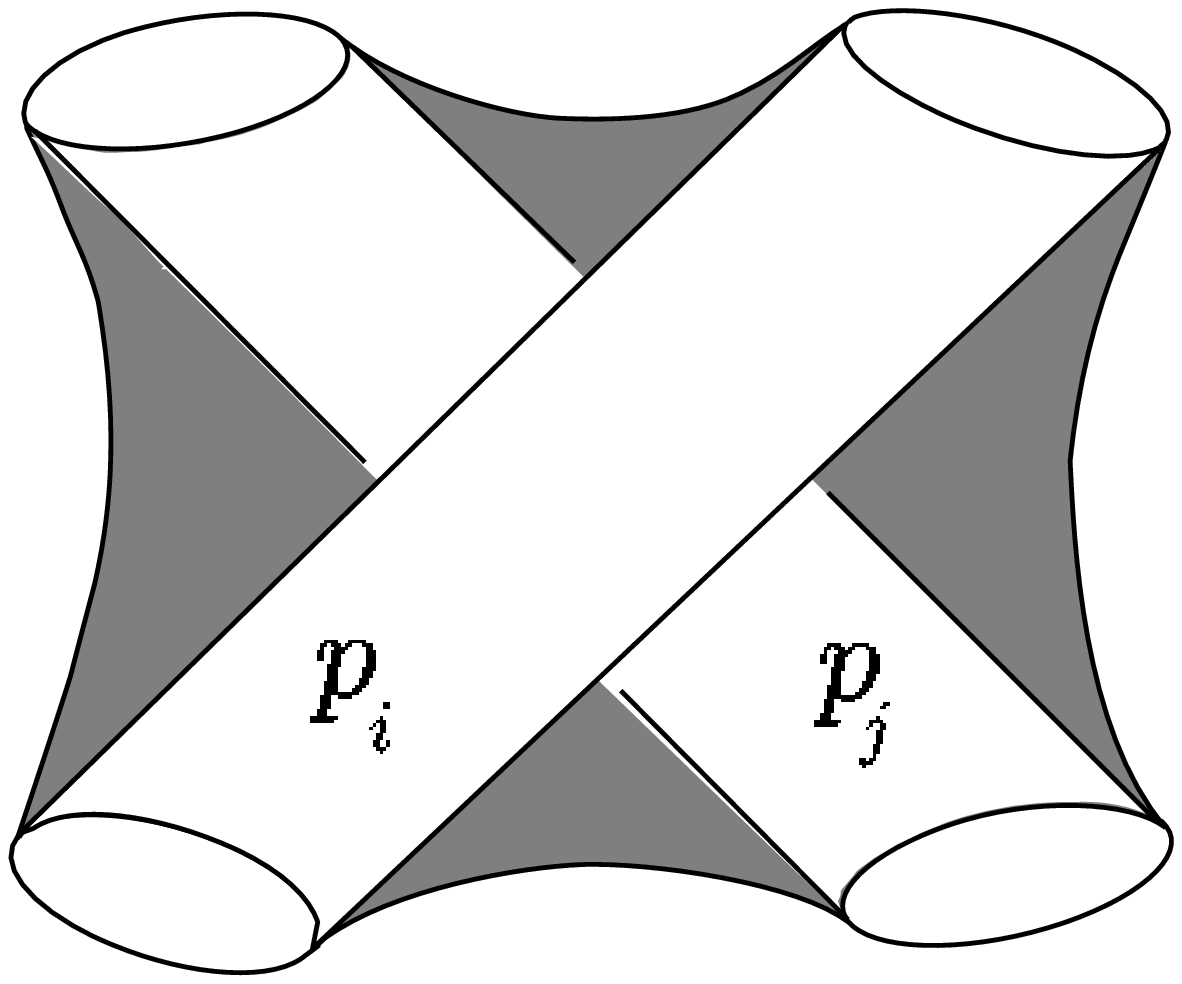}}
\put(13,1.9){$\Biggr )$}
\put(14,-0.6){\includegraphics[height=3.5cm]{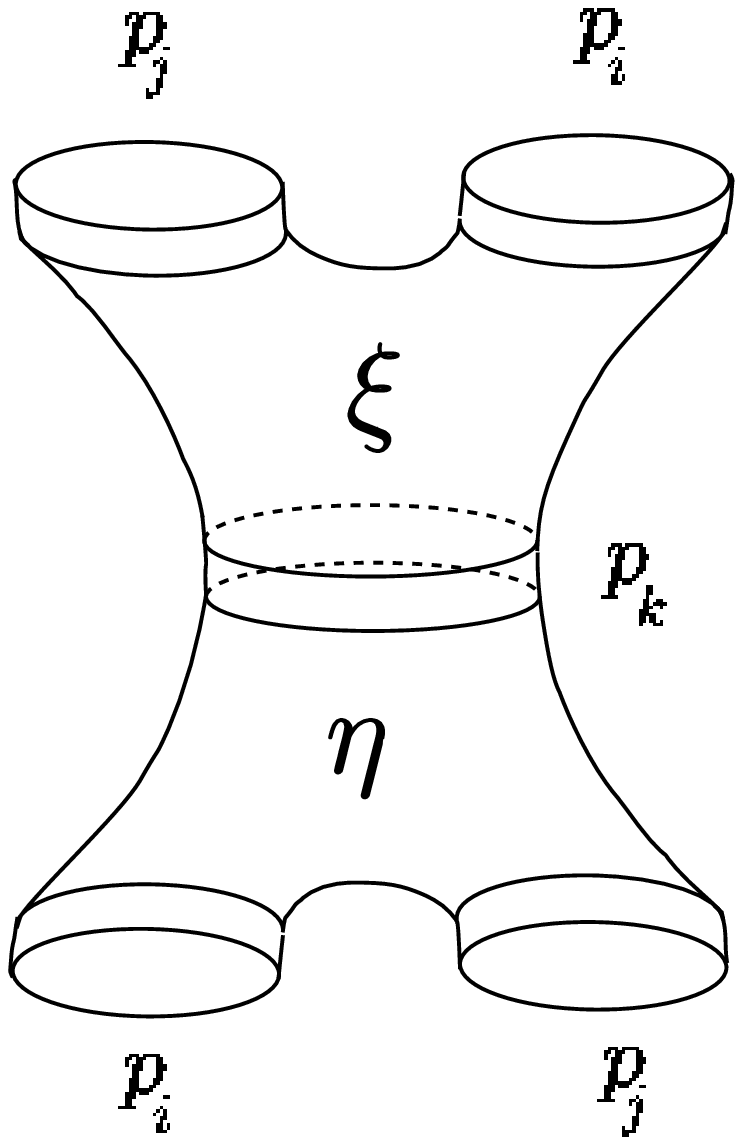}}
\end{picture}
$$

\end{defn}

The following theorem is due to Ocneanu. 
\begin{thm}
$Z_{\cal M}(S^1 \times S^1)$ is naturally identified with 
the center of Tube$({\cal M})$.
\end{thm}

Thanks to this theorem, one can consider the mapping class group 
of the torus $S^1 \times S^1$ as represented on the center 
of Tube$({\cal M})$. Then, we have the following theorem. 
\begin{thm}
$\{p_0, \cdots , p_d \}$ is a Verlinde basis of $Z(S^1 \times S^1)$. 
\end{thm}

Hence, for a closed oriented 3-manifold $V$ obtained from Dehn surgery 
along a link $L$, we have the Dehn surgery formula for the Turaev-Viro-Ocneanu 
TQFT:
$$
Z_{{\cal M}}(V)=\sum_{i_1,\cdots, i_m=0}^d S_{0i_1} \cdots S_{0i_m} 
J(L;i_1, \cdots, i_m)
$$
where $Z(S)p_i = \sum_{j=0}^d S_{ij}p_j$, ($i=0, \cdots ,d$). 

There are some applications of this theorem and we list them as 
corollaries. 

Note that, since the tube system is a modular category in the sense of 
Turaev, we can make a Reshetikhin-Turaev invariant of 3-manifolds out of it.

\begin{cor}
Let $V$ be a closed oriented 3-dimensional manifold. Then, 
Turaev-Viro-Ocneanu invariant of $V$ is equal to the 
Reshetikhin-Turaev invariant constructed from the tube system. 
\end{cor}

We remark that A. Ocneanu also has the above statement in \cite{O}. 

\begin{cor}
The Dijkgraaf-Witten invariant \cite{DW, W} of finite group $G$ with 
3-cocycle $\omega$ is equal to Altshuler-Coste invariant \cite{AC} of 
the quantum double $D^\omega(G)$. 
\end{cor}

\begin{cor}
If we further assume our $C^*$-tensor category $\cal M$ is modular, 
then $Z_{\cal M} (V)=\tau_{\cal M} (V) \overline{\tau_{\cal M} (V)}$, 
where $\tau_{\cal M} (V)$ is the Reshetikhin-Turaev invariant of 
a closed oriented 3-dimensional manifold $V$ constructed from 
$\cal M$ and the bar means the complex conjugation. 
\end{cor}

The last corollary follows because our tube system splits 
into ${\cal M} \otimes \overline{\cal M}$ under the assumption that 
${\cal M}$ is modular. 

\section{Computations of Turaev-Viro-Ocneanu invariants}
In the analysis of the Longo-Rehren subfactor, which corresponds to 
the center construction in the sense of Drinfel'd, M. Izumi 
formulated a tube algebra in terms of sectors for a finite closed 
system of endomorphisms of a type III subfactor in a slightly different 
way \cite{Izumi1}, although Izumi's tube algebra is isomorphic 
to Ocneanu's one as $C^*$-algebras. He also explicitly gave an action 
of $SL_2(\Bbb{Z})$ on the center of the tube algebra in the language 
of sectors, and derived some formulas about Turaev-Viro-Ocneanu invariants 
of lens spaces for concretely given subfactors. However, it is not clear how 
$S$ and $T$-matrices in Izumi's sector theory associated with the 
Longo-Rehren subfactor and $S$ and $T$-matrices in Turaev-Viro-Ocneanu 
$(2+1)$-dimensional TQFT are related. Actually, we can prove that for 
any subfactor an isomorphism between two tube algebras in Izumi's sector 
theory and Ocneanu's one, which is derived from Turaev-Viro-Ocneanu 
$(2+1)$-dimensional TQFT, induces a conjugate linear isomorphism as matrix 
representations of $SL_2(\Bbb{Z})$. Hence, we can compute a plenty of 
Turaev-Viro-Ocneanu invariants of 3-manifolds based on Izumi's data 
\cite{Izumi2} of $S$ and $T$-matrices. 

As a summary of the above argument, we state the following proposition. 

\begin{prop}
The representation of $SL_2(\Bbb{Z})$ constructed from the 
Turaev-Viro-Ocneanu TQFT are complex conjugate to ones introduced by Izumi. 
\end{prop}

We have the following lists by partially using the Maple software for 
analytical computation in mathematics. In the list, we use the notation 
$Z_{\cal M}(V)$ in such a way that it means that the Turaev-Viro-Ocneanu 
invariant of a 3-manifold $V$ constructed from the data $\cal M$, which is 
obtained from a subfactor $N \subset M$. (See also Section 3 for these 
notations.)

\par \bigskip 
\noindent \underline{the $D_5^{(1)}$-subfactor}
\par 
\par 
$Z_{\text{\normalsize $D_5^{(1)}$}}(L(3,1))=\overline{Z_{\text{\normalsize $D_5^{(1)}$}}(L(3,2))}=\dfrac{1+2w^2}{6}$, where $w^3=1$
\par \medskip 
\par 
$Z_{\text{\normalsize $D_5^{(1)}$}}(L(5,1))=Z_{\text{\normalsize $D_5^{(1)}$}}(L(5,2))=Z_{\text{\normalsize $D_5^{(1)}$}}(L(7,1))=Z_{\text{\normalsize $D_5^{(1)}$}}(L(7,2))=\dfrac{1}{6}$

\par \bigskip 
\noindent \underline{an $E_6$-subfactor} \cite{SuzukiWakui}
\par 
$Z_{\text{\normalsize $E_6$}}(L(p,1))$=$\dfrac{1}{12}\{ ((-1)^p+1)e^{-\frac{p\pi i}{3}}+2e^{-\frac{5p\pi i}{6}}+i^p+2(-1)^p+5\} $
\par \medskip 
$Z_{\text{\normalsize $E_6$}}(L(p,2))$=$\dfrac{1}{4}+\dfrac{(-1)^{\frac{p+1}{2}}i}{12}-\dfrac{\sqrt{3}+i}{12}e^{-\frac{(p+1)\pi i}{6}}$

\par 
\bigskip 
\begin{tabular}{c|cccc}
$(p,q,r)$ &  $(2,3,5)$ & $(2,3,7)$ & $(2,5,7)$ & $(3,5,7)$ \\ 
\hline 
$Z_{\text{\normalsize $E_6$}}(M(p,q,r))$ & $\frac{(6+2\sqrt{3})+(3-3\sqrt{3})i}{12}$ & $\frac{(6+2\sqrt{3})+(3-3\sqrt{3})i}{12}$ & $\frac{-\sqrt{3}+9+6i}{12}$ & $\frac{2-\sqrt{3}i}{2}$ \\ 
\end{tabular}

\par \bigskip 
\noindent \underline{the generalized $E_6$-subfactor with $G=\mathbb{Z}/3\mathbb{Z}$} --- see footnote\footnote{Values corrected 21 March 2003}
\par 
$Z_{\text{\normalsize $E_6,\mathbb{Z}/3\mathbb{Z}$}}(L(3,1))=\overline{Z_{\text{\normalsize $E_6,\mathbb{Z}/3\mathbb{Z}$}}(L(3,2))}=\dfrac{7-\sqrt{7}i}{14}$
\par \medskip 
$Z_{\text{\normalsize $E_6,\mathbb{Z}/3\mathbb{Z}$}}(L(5,1))=Z_{\text{\normalsize $E_6,\mathbb{Z}/3\mathbb{Z}$}}(L(5,2))=\dfrac{7-\sqrt{21}}{42}$
\par \medskip 
$Z_{\text{\normalsize $E_6,\mathbb{Z}/3\mathbb{Z}$}}(L(7,1))=Z_{\text{\normalsize $E_6,\mathbb{Z}/3\mathbb{Z}$}}(L(7,2))=\dfrac{1+\sqrt{3}i}{6}$

\par \bigskip 
\noindent \underline{the generalized $E_6$-subfactor with $G=\mathbb{Z}/4\mathbb{Z}$}
\par 
$Z_{\text{\normalsize $E_6,\mathbb{Z}/4\mathbb{Z}$}}(L(3,1))=Z_{\text{\normalsize $E_6,\mathbb{Z}/4\mathbb{Z}$}}(L(3,2))=\dfrac{2+\sqrt{2}}{16}$
\par \medskip 
$Z_{\text{\normalsize $E_6,\mathbb{Z}/4\mathbb{Z}$}}(L(5,1))=Z_{\text{\normalsize $E_6,\mathbb{Z}/4\mathbb{Z}$}}(L(5,2))=\dfrac{2+\sqrt{2}}{16}$
\par \medskip 
$Z_{\text{\normalsize $E_6,\mathbb{Z}/4\mathbb{Z}$}}(L(7,1))=Z_{\text{\normalsize $E_6,\mathbb{Z}/4\mathbb{Z}$}}(L(7,2))=\dfrac{2-\sqrt{2}}{16}$

\par \bigskip 
\noindent \underline{the generalized $E_6$-subfactor with $G=\mathbb{Z}/5\mathbb{Z}$} --- see footnote${}^1$
\par 
$Z_{\text{\normalsize $E_6,\mathbb{Z}/5\mathbb{Z}$}}(L(3,1))=Z_{\text{\normalsize $E_6,\mathbb{Z}/5\mathbb{Z}$}}(L(3,2))=\dfrac{1-\sqrt{5}}{10}$
\par \medskip 
$Z_{\text{\normalsize $E_6,\mathbb{Z}/5\mathbb{Z}$}}(L(5,1))=\dfrac{1}{3},\ Z_{\text{\normalsize $E_6,\mathbb{Z}/5\mathbb{Z}$}}(L(5,2))=\dfrac{2}{3}$
\par \medskip 
$Z_{\text{\normalsize $E_6,\mathbb{Z}/5\mathbb{Z}$}}(L(7,1))=Z_{\text{\normalsize $E_6,\mathbb{Z}/5\mathbb{Z}$}}(L(7,2))=\dfrac{3+\sqrt{5}}{30}$

\par 
\bigskip 
\noindent \underline{the generalized $E_6$-subfactor with $G=\mathbb{Z}/2\mathbb{Z}\times \mathbb{Z}/2\mathbb{Z}$}
\par 
$Z_{\text{\normalsize $E_6,\mathbb{Z}/2\mathbb{Z}\times \mathbb{Z}/2\mathbb{Z}$}}(L(3,1))=Z_{\text{\normalsize $E_6,\mathbb{Z}/2\mathbb{Z}\times \mathbb{Z}/2\mathbb{Z}$}}(L(3,2))=\dfrac{2+\sqrt{2}}{16}$
\par \medskip 
$Z_{\text{\normalsize $E_6,\mathbb{Z}/2\mathbb{Z}\times \mathbb{Z}/2\mathbb{Z}$}}(L(5,1))=Z_{\text{\normalsize $E_6,\mathbb{Z}/2\mathbb{Z}\times \mathbb{Z}/2\mathbb{Z}$}}(L(5,2))=\dfrac{2+\sqrt{2}}{16}$
\par \medskip 
$Z_{\text{\normalsize $E_6,\mathbb{Z}/2\mathbb{Z}\times \mathbb{Z}/2\mathbb{Z}$}}(L(7,1))=Z_{\text{\normalsize $E_6,\mathbb{Z}/2\mathbb{Z}\times \mathbb{Z}/2\mathbb{Z}$}}(L(7,2))=\dfrac{2-\sqrt{2}}{16}$
\par 

\par \bigskip 
\noindent \underline{the Haargerup subfactor}
\par 
$Z_{Haagerup}(L(7,1))=Z_{Haagerup}(L(7,2))=\dfrac{13+3\sqrt{13}}{78}$
\par \medskip 
$Z_{Haagerup}(M(2,3,5))=-\dfrac{\sqrt{13}}{26}+\dfrac{7}{6}$

%
%

\Addresses\recd
\end{document}